\input amstex
\documentstyle{amsppt}
\def\conj{\roman{conj}}
\def\Int{\mathop{\roman{Int}}\nolimits}
\def\Cl{\mathop{\roman{Cl}}\nolimits}
\def\nologo{\let\logo@\relax}
\let\ge\geqslant
\let\le\leqslant
\def\C{{\Bbb C}}
\def\R{{\Bbb R}}

\def\Q{{\Bbb Q}}

\def\Cp#1{\C\roman P^{#1}}
\def\Rp#1{\R\roman P^{#1}}
\def\QR{$Q$}
\let\tm\proclaim
\let\endtm\endproclaim

\def\s{\sigma}
\def\e{\varepsilon}
\def\lk{\mathop{l}\nolimits}
\def\d{\partial}
\def\G{\Gamma}
\def\om{\omega}
\def\oW{W^\circ}
\def\ra{\rangle}
\def\la{\langle}
\def\oo{\varnothing}

\magnification\magstep1
\NoBlackBoxes
\NoRunningHeads
\nologo

\topmatter
\title
Generalization of Arnold---Viro Inequalities
for real singular algebraic curves 
\endtitle
\author
Sergey Finashin
\endauthor
\address
Middle East Technical University, Ankara, Turkey
\endaddress{}
\email
serge\,\@\,rorqual.cc.metu.edu.tr, serge\,\@\,salamis.emu.edu.tr, 
finash\,\@\,pdmi.ras.ru
\endemail
\endtopmatter

\document
\heading
\S1. Introduction
\endheading
\subheading{1.1. Petrovskii and Arnold inequalities}
The subject of the present paper is related to the 16th Hilbert problem,
namely, to the question of possible distribution of ovals of plane
real algebraic curves. Recall that
the particular Hilbert's question of possibility of a
non-singular sextic having 11 ovals 
(the maximal possible number by Harnack theorem)
lying separately was answered negatively. The prohibition follows
from the Petrovskii inequality, which states that
$$\chi(W)\ge-\frac32k(k-1),$$ 
where
$
W=\{[x:y:z]\in\Rp2\,|\,f(x,y,z)\ge0\}
$
and $f$ is a real form of degree $2k$ defining a non-singular curve;
we denote by $A$ the complex point set of the latter
and by $A_\R$ its real part, $A\cap\Rp2$.
The change of the sign of $f$ substitutes $W$ for the complementary
region; thus we get the other Petrovskii inequality,
$\chi(W)\le\frac32k(k-1)+1$.

The modern approach for studying the topology of $A_\R$ in $\Rp2$,
developed in works of V.~A.~Rokhlin and V.~I.~Arnold,
is to use its close relationship with
the topology of the double plane, $\pi\:X\to\Cp2$,
branched along $A$.
One can define $X$ by the equation
$f(x,y,z)=t^2$ in a quasi-homogeneous complex projective $3$-space,
so it is a real algebraic surface;
we denote
by $\conj\: X\to X$ the involution of complex conjugation in
$X$ and by $X_\R$ the fixed point set of $\conj$, which is obviously
projected by $\pi$ into $W$.
Arnold \cite{A} noticed that the Petrovskii inequality can be obtained
by analyzing the negative inertia indices of the intersection form 
restricted to the eigenspaces of the involution
$\conj_*\:H_2(X;\R)\to H_2(X;\R)$.
This idea and its original proof \cite{A} can be interpreted 
in terms of the quotient manifold $Y=X/\conj$ as follows.
We subtract the Hirzebruch signature formula 
from the Riemann--Hurwitz formula applied to the branched coverings
$X\to Y$ and $X\to\Cp2$
$$
\align
2\chi(Y)-\chi(X_\R)&=\chi(X)=2\chi(\Cp2)-\chi(A)=4k^2-6k+6\\
2\s(Y)-X_\R\circ X_\R&=\s(X)=2\s(\Cp2)-A\circ A=2-2k^2,
\endalign
$$
($X_\R\circ X_\R$ and $A\circ A$ are self-intersection numbers in $X$)
and take into account that $b_1(X)=b_1(Y)=0$, and
$X_\R\circ X_\R=-\chi(X_\R)$ (the latter is due to the anti-isomorphism 
between the normal and the tangent bundle to $X_\R$). We obtain
$$
b_2^-(Y)-\frac12\chi(X_\R)=\frac12(b_2^-(X)-1)=\frac32k(k-1).
$$
Furthermore, $\chi(X_\R)=2\chi(W)$, and
$0\le b_2^-(Y)=\frac32k(k-1)+\chi(X_\R)$ can be seen as 
the Petrovskii inequality.

This point of view resulted about in a strengthening
of Petrovskii inequalities, namely
the Arnold inequalities, which can be read (after $b_2^\pm(Y)$ is 
similarly expressed) as
$$
\align
\s_\pm&\le b_2^\pm(Y)\\
\s_\pm+\s_0&\le b_2^\pm(Y)+1,
\endalign
$$
where $\s_+$, $\s_-$, $\s_0$ are the numbers of the connected
components of $X_\R$ with negative, positive
and zero Euler characteristic respectively.

\subheading{1.2. Generalizations of the Arnold inequalities for nodal
curves}
V.I.Zvonilov \cite{Z} extended the Arnold inequalities to non-singular 
curves $A$ of odd
degrees by applying a version of these inequalities to the curve obtained
from $A$ by adding a line.
This involved consideration of real curves with nodal singularities.
The extension of the Arnold inequalities to arbitrary plane real nodal
curves was formulated by O.Ya.Viro \cite{V1}. In the case of nodal curves,
the double plane $X$ is no longer
a manifold, but is a rational homology manifold, so the approach of
Arnold still works. However, the numbers of components, $\s_\pm$, $\s_0$,
should be replaced by the inertia indices $\s_\pm(q)$, $\s_0(q)$,
of the form induced on $H_2(X_\R)$ from the intersection
form in $X$. 
Viro described this form in combinatorial terms
and determined the upper estimates
in the generalized Arnold inequalities.
These estimates were further improved by V.\,M.\,Kharlamov and O.\,Ya.\,Viro.
A modification of the inequalities treating nodal surfaces was also
formulated by Kharlamov \cite{Kh}.
 The proof of Viro \cite{V2} was not
published, but a slightly simplified version of the proof,
in which the original constructions are descended into the quotient $Y$,
can be found in \cite{F1}.

\subheading{1.3. The results}
In this paper, I present a generalization of the Arnold--Viro inequalities
to the case of a curve with the other isolated singularities with which
 $X$ is still a rational homology manifolds.
The detailed proofs, some applications and
further generalizations and refinements (for curves 
with more general singularities on non-singular 
surfaces and for singular surfaces,
estimates of the sharpness of the inequalities)
 are to appear in \cite{F2}.

\heading
\S2. 
Generalized Arnold--Viro inequalities
\endheading
\subheading{2.1. Linking forms}
Assume that $U$ is a sufficiently
small compact regular, thus cone-like, neighborhood of an
isolated surface singularity at point $p$; we call it
\QR-singularity if $M=\d U$ is a rational homology sphere.
It is well known (and trivial) that a singularity is \QR-singularity
if and only if its resolution graph is a tree with all vertices
represented by spheres. For singularities of functions
\QR-property is equivalent to non-degeneracy of the Milnor form.

Assume now that the singularity at $p$ is real and $U$ is invariant
under the complex conjugation, $\conj$.
Denote by $U_\R$ the fixed point set of $\conj$ in $U$.
We obtain a link $L=\d U_\R=U_\R\cap M$ in $M$, which has a natural framing
given by the outward normal vector field to $L$ in $U_\R$
multiplied by $i$.
Thus we obtain {\it the linking form of a real \QR-singularity},
$\lk_p\: H_1(L)\to\Q$.
Denote by $\lk_p+\roman I$ the form obtained from $\lk_p$ by adding
$1$ to the self-linking numbers
of the components of $L$ and preserving the linking coefficients
between distinct components.

\subheading{2.2. The forms $q_p$ for real curve singularities}
Consider now a curve $A\subset\C^2$ defined by a polynomial $f(x,y)$,
with real coefficients. Its  singularity
at $p=(x_0,y_0)\in\C^2$ will be called \QR-singularity if its suspension
(that is the corresponding singularity at $\hat p=(x_0,y_0,0)$
on the surface $\{f(x,y)=z^2\}\subset\C^3$) is  \QR-singularity.
Assume now that $p\in\R^2$ and
consider a neighborhood $U\ni\hat p$, 
$U=\{(x,y,z)\in\C^2\,|\,|x-x_0|^2+|y-y_0|^2\le\e,|z|^2\le\e_1\}$,
for $0<\e<\!\!\!<\e_1<\!\!\!<1$.
The link $L$ defined above is projected by $\pi\:\C^3\to\C^2$,
$(x,y,z)\mapsto(x,y)$, into $S_p^+=S_p\cap\{(x,y)\in\R^2\,|\,f(x,y)\ge0\}$,
where $S_p$ is the circle of radius $\e$ around $p$.
Fix an orientation of $S_p$; then
split $L$ into halves, $L^\pm=\{(x,y,z)\in L\,|\, \pm z\ge0\}$,
and orient $L$ so that $\pi$ preserves (reverses) the orientation
of $L^+$ ($L^-$).
Denote by $q_p\:\Cal H_p\to\Q$ the quadratic form on 
$\Cal H_p=H^0(S_p^+)$ obtained as the pull back of $\lk_p+\roman I$
via the product of the homomorphism
$\pi^*\:H^0(S_p^+)\to H^0(L)$ and
the Poincare duality $H^0(L)\to H_1(L)$.
Note that $q_p$ is independent of
the orientation of $S_p$.
If $p\notin\Cl(\Int W)$
(i.e., $\rho=0$ and $f$ is negative around $p$), then 
$S_p^+=\oo$, and we put $\Cal H_p=0$, $q_p=0$.

Assume that a simple closed curve $\om\subset A_\R$ contains $p$;
then $S_p$ is split by $\om\cap S_p$ in two arcs.
Fix one of them and change the above orientation of $L$
on the components which are projected
into this arc. The new fundamental class
$[L]\in H_1(L)$ changes the duality $H^0(L)\to H_1(L)$, so the pull back
of $\lk_p+\roman I$ is another
quadratic form, $q_{p,\om}\:\Cal H_p\to\Q$,
which is independent of the choice of the arc and the
orientation of $S_p$.

\subheading{2.3. Partition forms}
Consider now 
a real curve $A\subset\Cp2$ of degree $2k$ defined
by a homogeneous polynomial $f(x,y,z)$.
For odd $k$ denote by $W_1,\dots,W_n$ the partition
components of $W=\{[x:y:z]\in\Rp2\,|\, f(x,y,z)\ge0\}$
(the closures of the connected components of
the interior, $\Int(W)$).
For even $k$ denote by $W_1,\dots,W_n$ the orientable partition
components of $W$. 
Put $W^\circ=\cup_{i=1}^n\Int W_i$.
If $k$ is even, then choose a simple closed piecewise smooth 
non-contractible in $\Rp2$ curve $\om\subset\Rp2-\oW$.
Such a curve can be chosen inside the non-orientable component of
$W$ if such a component exists; otherwise the inclusion homomorphism
$H_1(A_\R)\to H_1(\Rp2)$ is non-trivial and we can choose
$\om\subset A_\R$.
For odd $k$ put $\om=\oo$.
Denote by $S$ the set of all singularities of $A$ and put
$S_\R=S\cap\Rp2$.

Assume for the rest of this section
that all the real singular points of $A$, $p\in S_\R$, are \QR-singular.
Denote by $w^i\in H^0(\oW)$ \ ($w^i_p\in \Cal H_p$), $i=1,\dots, n$,
 the classes represented by
the cochains equal to $1$ on $\Int W_i$ \ ($W_i\cap S_p$)
and to $0$ on the rest. Note that $w^i_p=0$ if $p\notin W_i$.

We define the {\it partition form} $q\:H^0(\oW)\to\Q$ as follows
$$
\align
q(w^i,w^j)&=\sum_{p\in S_\R-\om}q_p(w^i_p,w^j_p)+
\sum_{p\in S_\R\cap\om}q_{p,\om}(w^i_p,w^j_p),
\phantom{AAAAAAAA}
\text{if }\ 1\le i,j\le n, i\ne j\\
q(w^i,w^i)&=\sum_{p\in S_\R-\om}q_p(w^i_p,w^i_p)+
\sum_{p\in S_\R\cap\om}q_{p,\om}(w^i_p,w^i_p)-2\chi(\Int W_i),\ 
\text{if }\ 1\le i\le n.
\endalign
$$

\subheading{2.4. Generalized Arnold--Viro inequalities}
Let $r$ denote the number of irreducible components of $A$,
$\nu=0$ if all irreducible components of $A$ have even degree and 
$\nu=1$ otherwise;
$\frak m(p)$ the Milnor form of the singularity at $p\in S$,
and $\mu_\pm=\s_\pm(\frak m)$,
where $\frak m=\oplus_{p\in S}(\frak m(p))$ is the total Milnor form of $A$,
$\s_\pm$ denote the inertia indices of quadratic forms and
$\s_0$ denotes the nullity.

\tm{Theorem A}
(Generalized Arnold--Viro inequalities)
$$
\align
\s_+(q)&\le
\frac12(k-1)(k-2)-\frac12\mu_+\\
\s_+(q)+\s_0(q)&\le 
\frac12(k-1)(k-2)-\frac12\mu_+
+(r-\nu)\\
\s_-(q)&\le
\frac32k(k-1)+\frac12\chi(X_\R)-\frac12\mu_-\\
\s_-(q)+\s_0(q)&\le 
\frac32k(k-1)+\frac12\chi(X_\R)-\frac12\mu_-+(r-\nu)\\
\endalign
$$
\endtm

\subheading{2.5. Scheme of the proof}
By Edmonds theorem 
the fixed point set
of an involution on a spin manifold carries a naturally defined
semi-orientation
(a pair of the opposite orientations),
provided the involution preserves the orientation and the spin structure.
For odd $k$ the complement
$X-S_X$ of the singularity, $S_X\subset X$,
is a spin simply connected manifold, thus,
$X_\R-S_X$ obtains a semi-orientation. 
If $k$ is even, then we can define the relative with respect to $\om$
spin semi-orientation of the union of the
components $\G_i=\pi^{-1}(W_i)\subset X_\R$, $i=1,\dots,n$,
(outside the singularity)
cf. \cite{F1}.
Let us fix any of the two orientations provided by
the spin (relative spin) semi-orientation
and denote by $[\G_i]$ the fundamental classes of $\G_i$.
Note that the quotient $Y$ is a rational homology manifold 
since $X$ is.

Denote by
 $\la\,,\,\ra_X$ and $\la\,,\,\ra_Y$ the intersection forms in $X$ and $Y$.

The following lemmas imply Theorem A.

\tm{Lemma 1} For any $1\le i,j\le n$
$$
q(w^i,w^j)=\la[\G_i],[\G_j]\ra_X=\frac12\la[\G_i],[\G_j]\ra_Y
$$
\endtm

\tm{Lemma 2} 
$$
\align
b_2^+(Y)&=\frac12(k-1)(k-2)-\frac12\mu_+\\
b_2^-(Y)&=\frac32k(k-1)+\frac12\chi(X_\R)-
\frac12\mu_-
\endalign
$$
\endtm

\tm{Lemma 3}
$\dim\ker(H_2(X_\R;\R)\to H_2(Y;\R))\le r-\nu$
\endtm

The versions of the above
lemmas for non-singular and for nodal curves $A$ are well known.
To prove Lemma 3 in our more general setting we do not need any changes
in the arguments, which can be found for instance in \cite{F1}.
Proofs of Lemmas 1---2 are straightforward calculations, although
rather less trivial then for nodal curves.

\heading
{\S3. Computation of the forms $q_p$}
\endheading
\subheading{3.1. The local partition forms}
Assume that a real polynomial $f(x,y)$ defines
an isolated singularity at $p\in\R^2$; denote by
$\mu$ its Milnor number and by $\rho$ 
 the number of real branches at $p$.
Consider a small real deformation (morsification) $\widetilde f$, of $f$
such that the singularity at $p$ produces the maximal possible number,
$\delta=\frac12(\mu+\rho-1)$, of real hyperbolic nodes on the curve 
$\{\widetilde f=0\}$, cf. \cite{M},
which all lie inside a closed $\e$-disk
$B_p\subset\R^2$, $0<\e<\!\!\!<1$, around $p$
(recall that a hyperbolic node is a singularity of type  $A_1$,
whose branches are real).
The method for constructing such a deformation can be found in
\cite{AC}, \cite{GZ}.
We put $\widetilde W_p=\{\widetilde f(x,y)\ge0\}\cap B_p$,
and denote by 
$\widetilde W_1,\dots,\widetilde W_{m}$ the closures of 
the connected components of $\Int(\widetilde W)$.
Let $\widetilde{\Cal H}_p=H^0(\Int(\widetilde W))$.
Assume that $\widetilde W_1,\dots, \widetilde W_n$ lie in the interior
of $B_p$ and $\widetilde W_{n+1},\dots, \widetilde W_m$
intersect $\d B_p$. 
If $\rho\ge1$, then $m=n+\rho$. If $\rho=0$, then $m=n+1$
if $f$ is positive around $p$ and  $m=n$ if negative.
To define the following quadratic form,
$\widetilde q\:\widetilde{\Cal H}_p\to\Q$, we use the local versions
of the formulas in 2.3, which are apparently
the local versions of Viro formulas \cite{V1},
since $\widetilde f$ has no other singularities except hyperbolic nodes.
Namely,
$$
\align
\widetilde q(w^i,w^j)&=\frac12
\#(\widetilde W_i\cap \widetilde W_j);
\phantom{AAAAAAAAAAA}\,\,
\text{for $1\le i,j\le m$, $i\ne j$}\\
\widetilde q(w^i,w^i)&=
\frac12\#(\widetilde W_i\cap\Cl(\widetilde W-\widetilde W_i))
-2\chi(\widetilde W_i)
\text{\ \ for $1\le i\le m$}\\
\endalign
$$
where $\#$ denotes the number of points and
$w^i\in \widetilde{\Cal H}_p$ are the generators representing
$\Int(\widetilde W_i)$.

Denote by $E$ the subspace of $\widetilde{\Cal H}_p$
generated by $w^1,\dots,w^n$.
Note that $\Cal H_p$ introduced in 2.2 
can be naturally identified with the subspace of  $\widetilde{\Cal H}_p$ 
generated by $w^{n+1},\dots,w^m$.

\tm{Theorem B}
Assume that we have a \QR-singularity at $p$. Then 
\roster
\item
the restriction of $\widetilde q$ to $E$
is non-degenerated; in particular, we have a decomposition
$\widetilde{\Cal H}_p=E\oplus E^\bot$, where $E^\bot$ is the orthogonal
complement to $E$ with respect to $\widetilde q$,
and the orthogonal projection
$pr\:{\Cal H}_p\to E^\bot$ is an isomorphism;
\item
if  $\rho\ge1$  or $\rho=0$ and $f$ is positive around $p$, then
the form $q_p$ introduced in 2.2 is equal to 
$\widetilde q\circ pr+2{\roman I}$
in $\Cal H_p\subset\widetilde{\Cal H}_p$,
where ${\roman I}$ denotes the quadratic form defined by the identity matrix
in the basis $w^{n+1},\dots,w^m$
(if  $\rho=0$ and $f$ is negative around $p$,
then obviously $\Cal H_p=E^\bot=0$ and $q_p=\widetilde q\circ pr=0$).
\endroster
\endtm

Note that this theorem really gives an explicit method to
calculate the forms $q_p$, because $\widetilde q\circ pr$
is obtained from $\widetilde q$
by the usual algorithm of completing the squares of
$w^1,\dots,w^n$ and taking the rest, which is 
a form in variables $w^{n+1},\dots,w^m$.

The forms $q_{p,\om}$ are obviously determined by $q_p$.

\subheading{3.2. The case of simple singularities}
As an illustration
we present below the matrices $M=M(q_p)$ of the forms $q_p$ for
the simple real surface singularities, which can be easily computed
applying the above algorithm to the well known (cf. \cite{AC}, \cite{GZ})
real morsifications of these singularities.

Note that we may not specify the correspondence between 
variables, $w^i$, and the entries of matrices due to the symmetry
in all the cases except
tree: for $D_{2n+2}$, $D_{2n+3}$ and $E_7$.
The correspondence in the exceptional cases is determined by the rule:
the lesser diagonal entry, $q_p(w^i,w^i)$, corresponds to
the angle $0$ between the real branches
at $(0,0)$ which embrace the corresponding region $W_i$.

$$
\align
A_{2n-1}:\phantom{A}\,
& f(x,y)=\cases
-x^{2n}+y^2
&\phantom{AA}\,\,\,
M=\phantom{AAA\,}
\pmatrix 
\frac n2&\frac n2\\ 
\frac n2&\frac n2
\endpmatrix
\\
x^{2n}-y^2;
&\phantom{AA}\,\,\,
M=\phantom{a}\,
\pmatrix 
\frac{2n-1}{2n}&\frac1{2n}\\ 
\frac1{2n}&\frac{2n-1}{2n}
\endpmatrix
\\
x^{2n}+y^2;
&\phantom{AA}\,\,\,
M=\phantom{AAAAA}\!\!\!
\pmatrix2n\endpmatrix
\\
-x^{2n}-y^2;
&\phantom{AA}\,\,\,
M=\phantom{AAAAA}\!
\pmatrix0\endpmatrix
\endcases
\\
A_{2n}:\phantom{AA}
 &f(x,y)=
\cases
\pm x^{2n+1}+y^2;
&\phantom{A}\,\,
M=\phantom{AAAA}\,
\pmatrix 2n \endpmatrix
\\
\pm x^{2n+1}-y^2;
&\phantom{A}\,\,
M=\phantom{AAA}\,\,
\pmatrix\frac{2n}{2n+1}\endpmatrix
\\
\endcases
\\
D_{2n+2}:\ \ &f(x,y)=
\cases
\pm x(x^{2n}-y^2);
&\phantom{A}
M=\pmatrix
1 &\frac12 &\frac12\\
\frac12 &\frac{n+1}2 &\frac n2\\
\frac12 &\frac n2 &\frac{n+1}2
\endpmatrix
\\
\pm x(x^{2n}+y^2);
&\phantom{A}
M=\phantom{AAA}\,
\pmatrix 4n-2\endpmatrix
\\
\endcases
\\
D_{2n+3}:\ \ &f(x,y)=
\cases
x(x^{2n+1}\pm y^2);
&M=\phantom{A}
\pmatrix
2n+1 &1 \\
1 &1
\endpmatrix
\\
-x(x^{2n+1}\pm y^2);
&M=\phantom{A}\!
\pmatrix
\frac{2n+3}4&\frac{2n+1}4 \\
\frac{2n+1}4 & \frac{2n+3}4
\endpmatrix
\endcases
\\
E_6:\phantom{AAA}
 &f(x,y)=
\cases
x^4\pm y^3;
&\phantom{AAAA}
M=\phantom{AAAAA}\!
\pmatrix6\endpmatrix
\\
-x^4\pm y^3;
&\phantom{AAAA}
M=\phantom{AAAAA}\!
\pmatrix2\endpmatrix
\\
\endcases
\\
E_7:\phantom{AAA} 
&f(x,y)=\ \ 
\pm y(x^3\pm y^2);
\phantom{AAAA}\!
M=\phantom{AAA}
\pmatrix
\frac72 &\frac32 \\
\frac32 &\frac32
\endpmatrix
\\
E_8:\phantom{AAA} &f(x,y)=\ \ 
\pm x^5\pm y^3;
\phantom{AAAAA}\,\,\,
M=\phantom{AAAAA}\!
\pmatrix8\endpmatrix
\endalign
$$


\Refs
\widestnumber\key{DFM}


\ref\key A
\by V. I. Arnold
\paper On the arrangement of ovals of real plane algebraic curves,
involutions on $4$-dimensional smooth manifolds and the arithmetic
of integral quadratic forms
\jour  Funct. Anal. Appl.
\vol 5
\issue 3
\yr 1971
\pages 169--178
\endref


\ref\key AC
\by N. A'Campo
\paper Le groupe de monodromie du d\'eploiement des singularit\'es
isol\'ees de courbes planes I
\jour Math. Ann.
\vol 213
\yr 1975
\pages 1--32
\endref


\ref\key F1
\by S. Finashin
\paper On topology of real plane algebraic nodal curves
\jour Algebra i Analiz (Transl. in St. Petersburg Math. J.)
\vol 8
\issue 6
\yr 1996
\pages 186--203
\endref


\ref\key F2
\by S. Finashin
\paper An invariant of isolated real surface singularities 
and generalization of Arnold--Viro inequalities
\jour to appear
\endref



\ref\key GZ
\by S. M. Gusein-Zade
\paper Intersection matrices for some singularities of functions
of two variables
\jour Funkct. Anal. Appl.
\vol 8
\issue 1
\yr 1974
\pages 11--15
\endref


\ref\key Kh
\by V. M. Kharlamov
\paper Topologia dejstvitel'nyh algebraicheskih mnogoobrazij,
(Topology of real algebraic varieties)
\inbook a survey in
I. G. Petrovskii, Izbrannyie trudyi, Systemyi Uravnenii
s Chastnyimi Proizvodnyimi, Algebraicheskaya Geometriya
(Collected works, Systems of Partial Differential Equations,
Algebraic Geometry)
\yr 1986
\bookinfo Moscow, Nauka
\pages 465--493
\endref

\ref\key V1
\by O. Ya. Viro
\paper
Obobschenie neravenstv Petrovskogo i Arnol'da na krivye s osobennostjami
(Generalization of the Petrovskii and Arnold
inequalities for curves with singularities)
\jour Uspechi mat. nauk
(Russian Math. Surveys)
\issue 3
\vol 33
\yr 1978
\pages  145--146
\endref


\ref\key  V2
\by O. Ya. Viro
\paper Handwritten notes
\jour
\vol
\yr 1978
\pages
\endref


\ref\key Z
\by V. I. Zvonilov
\paper Refinement of the Petrovskii--Arnold inequality
for curves of odd degree
\jour Funct. Anal. Appl.
\vol 13
\yr 1979
\pages 262--268
\endref

\endRefs
\enddocument